\documentclass[draft]{article}
\usepackage[english]{babel}
\usepackage{amssymb}
\usepackage{amsfonts}
\usepackage{amsmath}
\usepackage{amsthm}

\let\al=\alpha

\let\ga=\gamma
\let\om=\omega

\let\Ga=\Gamma

\let\wtd=\widetilde

\newcommand{\bbM}{\mathord{\mathbb M}}
\newcommand{\bbN}{\mathord{\mathbb N}}
\newcommand{\bbR}{\mathord{\mathbb R}}

\newcommand{\bbP}{\mathord{\mathbb P}}
\newcommand{\bbZ}{\mathord{\mathbb Z}}

\newtheorem{example}{Example}

\def\cP{\mathord{\cal P}}

\def\ze{\zeta}

\def\ga{\gamma}
\def\si{\sigma}
\def\la{\lambda}
\newcommand{\ka}{\kappa}

\let\td=\tilde
\let\wh=\widehat

\newcommand{\sspan}{\mathop{\rm span}\nolimits}
\newcommand{\card}{\mathop{\rm card}\nolimits}

\newcommand{\Aff}{\mathop{\rm Aff}\nolimits}
\newcommand{\codim}{\mathop{\rm codim}\nolimits}

\newcommand{\SO}{\mathop{\rm SO}\nolimits}
\newcommand{\OO}{\mathop{\rm O}\nolimits}

\newcommand{\GL}{\mathop{\rm GL}\nolimits}

\newcommand{\Tr}{\mathop{\rm Tr}\nolimits}
\newcommand{\LL}{\mathop{\rm L}\nolimits}

\newcommand{\Int}{\mathop{\rm Int}\nolimits}

\newtheorem{lemma}{Lemma}
\newtheorem{theorem}{Theorem}
\newtheorem{proposition}{Proposition}
\newtheorem{corollary}{Corollary}

\def\Int{\mathop{\hbox{\rm Int}}\nolimits}

\def\hskl{\vskip 0.35\baselineskip}
\def\tskl{\vskip 0.2\baselineskip}

\newcount\cdef
\def\nxtdef{\global\advance\cdef by 1\relax}
\def\remark#1{{\sc\hskl\noindent Замечание.\kern10pt}{\small #1}\tskl\noindent}

\def\clos{\mathop{\hbox{\rm clos}}\nolimits}

\title{\bf On flat complete causal Lorentzian manifolds}
\author{V.M.~Gichev, O.S.~Morozov}
\date{}
\begin{document}
\maketitle
\begin{abstract}
We describe up to finite coverings causal flat affine complete
Lorentzian manifolds such that the past and the future of any
point are closed near this point. We say that these manifolds are
strictly causal. In particular, we prove that their fundamental
groups are virtually abelian. In dimension 4, there is only one,
up to a scaling factor, strictly causal manifold which is not
globally hyperbolic. For a generic point of this manifold, either
the past or the future is not closed and contains a lightlike
straight line.
\end{abstract}

\section{Introduction}

Let $\bbM_n$ be Minkowski spacetime of dimension $n$ and $\cP_n$
be Poincare group of its affine automorphisms. We consider
manifolds that can be realized as quotient spaces $\bbM_n/\Ga$,
where $\Ga\cong\pi_1(M)$ is a discrete subgroup of $\cP_n$ whose
action on $\bbM_n$ is free and proper. These manifolds may be
defined in differential geometric terms as geodesically complete
Lorentzian manifolds with vanishing torsion and curvature. They
are characterized by the existence of an atlas of coordinate
charts with coordinate transformations in $\cP_n$, and the
completeness. The latter means that any affine mapping of a
segment in $\bbR$ to $M$ extends affinely to $\bbR$. Also, these
manifolds are complete affine manifolds with a compatible
Lorentzian metric.

Complete affine manifolds are studied since 60s. The following
question is known as Auslander's conjecture: is it true that
$\pi_1(M)$ for a compact flat complete affine $M$ is virtually
solvable? It remains unanswered but in many cases has the
affirmative answer (see \cite{Abels},  \cite{CDGM03} for details
and further references). If $M$ is not compact then $\pi_1(M)$ may
be free non-Abelian according to the remarkable example by
Margulis \cite{Mar83} that gives the negative answer to the
question of Milnor \cite{Milnor}. In paper \cite{Fried} by Fried,
Auslander's conjecture was proved for Lorentzian compact
4-manifolds; also, \cite{Fried} contains a description of causal
two ended quotient spaces $H/\Ga$, where $H$ is a subgroup of
$\cP_n$ that is simply transitive on $\bbM_n$.

In what follows, we assume that manifolds are Lorentzian, flat and
complete if the contrary is not stated explicitly. The Lorentzian
metric defines in each tangent space the pair of closed convex
round cones. Choosing one of them, we get locally a cone field. It
can be extended to the global cone field on $M$ or on a two-sheet
covering space of $M$. If $M$ admits no closed timelike curves
then $M$ is said to be causal. We describe up to finite coverings
complete flat causal Lorentzian manifolds which satisfy the
following additional condition: the past and the future of any
point $p\in M$ are closed near $p$. We say that these manifolds
are {\rm strictly causal\/}. Generic causal manifolds of the paper
\cite{Fried} are not strictly causal. On the other hand,
nontrivial strictly causal manifolds are never globally
hyperbolic. Manifolds of the latter class are well understood (the
paper \cite{Barbot} by Barbot contains many useful information on
them, including a classification, without the assumption of the
completeness). We give an explicit (parametric) description of
strictly causal manifolds up to finite coverings. The simplest
nontrivial example has dimension 4: $M=\bbM_4/\Ga$, where
$\Ga\cong\bbZ$. Its causal properties are somewhat surprising. The
manifold $M$ is the disjoint union
$$M=M^+\cup M^0\cup M^-,$$
where $M^+$ and $M^-$ are open and $M^0$ is closed (moreover,
$M^0$ is an affine hyperplane). For any $p\in M^+$, its future
$F_p$ is closed;  the past $P_p$ is not closed and contains
lightlike straight lines. Furthermore, $M$ admits an involutive
automorphism that transposes $M^+$ and $M^-$ (hence the future of
a point in $M^-$ is not closed). For $p\in M^0$, both $F_p$ and
$P_p$ are closed and contains no lightlike line. Also, any
strictly causal 4-manifold that is not globally hyperbolic is
homothetic to the manifold of this example.

In contrast to most of the cases above, the fundamental group
$\pi_1(M)$ of a strictly causal manifold $M$ is virtually abelian.
Moreover, $M$ can be finitely covered by the product of a torus
and Euclidean space. A finite cover of $M$ can be realized as a
(topologically trivial) vector bundle over $\bbM_k/\Ga$, where
$k\leq n$, $\bbM_k$ is embedded to $M_n$ as an affine
$\Ga$-invariant subspace, $\Ga$ is unipotent in $\bbM_k$, and the
holonomy group is defined by a bounded linear representation of
$\Ga$ in the fibre. Thus the problem is reduced to the case of
unipotent $\Ga$. Up to finite coverings, there are two types of
these manifolds. The first, elliptic, consists of manifolds that
admit a Riemannian metric such that the identical mapping is
affine. In other words, linear parts of transformations in $\Ga$
keep some positive definite quadratic form; then $\Ga$ contains a
finite index subgroup of translations by vectors in some spacelike
vector subspace. For manifolds of the second (parabolic) type,
$\Ga$ is a uniform lattice in a vector group $T$ whose action in
$\bbM_n$ is free and quadratic on $T$.

Some results in this article overlap with the recent paper
\cite{Barbot}).

\section{Preliminaries and statement of results}
Fixing the origin $o\in\bbM_n$, we identify $\bbM_n$ with the real
vector space $V$ equipped with a Lorentzian form $\ell$ of the
signature $(+,-,\dots,-)$. The set $\ell(v,v)\geq0$ is the union
of two closed convex round cones in $V$. Let $C$ be one of them.
The group $\Ga$ is assumed to act freely and properly in $V$ by
affine transformations whose linear parts keep $\ell$ and $C$. We
denote by $\ka$ the quotient mapping $\bbM_n\to M=V/\Ga$ and
define the {\it past} $P_p$ and the {\it future} $F_p$ of $p\in M$
by
\begin{eqnarray*}
P_p=\ka(v-C),\qquad F_p=\ka(v+C),\qquad v\in\ka^{-1}(p).
\end{eqnarray*}
Clearly, $P_p$ and $F_p$ do not depend on the choice of $v$. On
$M$, this defines the field of pointed convex cones
$C_p=d_v\ka(v+C)$, $v\in V$. The manifold $M$ is said to be {\it
causal} if $M$ admits no closed piecewise smooth timelike paths. A
smooth path $\eta$ is called {\it timelike} if $\eta'(t)\in
C_{\eta(t)}$ for all $t$; for {\it lightlike} paths, $\eta'(t)\in
\partial C_{\eta(t)}$ (note that lightlike paths are timelike).
The definition naturally extends to piecewise smooth paths.
Clearly, any timelike curve in $M$ can be lifted to a timelike
curve in $V$ and the projection of a timelike curve in $V$ is
timelike. The following observation is useful: $M$ is causal if
and only if
\begin{equation}\label{causa}
v\in V,\ \ga\in\Ga,\ \ga(v)\in v+C\quad \Longrightarrow\quad \ga=e,
\end{equation}
where $e$ denotes the identity of $\Ga$. Indeed, if $\ga\neq e$
then $\ga(v)\neq v$ and the projection to $M$ of the segment with
endpoints $v$ and $\ga(v)$ is a nontrivial closed timelike curve.
Conversely, we get a timelike curve $\td\eta$ in $V$ lifting a
timelike curve $\eta$ in $M$; hence $\td\eta$ lies in $v+C$ if it
starts at $v\in\ka^{-1}(p)$, $p\in M$. If $\eta$ is closed and
nontrivial then its endpoint is $\ga(v)$ for some
$\ga\in\Ga\setminus\{e\}$.

We say that $M$ is {\it strictly causal} if $M$ is causal and for
each $p\in M$ there exists a neighbourhood $U$ of $p$ such that
$U\cap F_p$ and $U\cap P_p$ are closed in $U$. For $\ga\in\Ga$,
set
\begin{eqnarray}
&\ga(v)=\la(\ga)v+\tau(\ga), \quad\mbox{where}\quad
 \la(\ga)\in\OO(\ell),\ \tau(\ga)\in V;\label{affde}\\
&G=\la(\Ga).\nonumber
\end{eqnarray}
Clearly, $\la:\,\Ga\to\OO(\ell)$ is a homomorphism and for all $\ga,\nu\in\Ga$
\begin{eqnarray*}
&\tau(\ga\circ\nu)=\la(\ga)\tau(\nu)+\tau(\ga).
\end{eqnarray*}
A linear subspace $X\subset V$ is called {\it spacelike} if $X\cap
C=\{0\}$; $X$ is spacelike if and only if $\ell$ is negative
definite on $X$. We say that $M,\Ga$ and $G$ are {\it unipotent}
if $G$ consists of unipotent linear transformations. The
classification problem can be reduced to the unipotent case.
\begin{theorem}\label{split}
A strictly causal flat complete Lorentzian manifold $M$ can be finitely covered by
$V/\Ga$, where $\Ga$ is abelian and satisfies following conditions:
there exists $\td o\in V$ and linear subspaces
$V_0,V_1\subseteq V$ such that $V=V_0\oplus V_1$ and
\begin{itemize}
\item[1)] the decomposition is $\ell$-orthogonal and $G$-invariant,
$\ell$ is Lorentzian in $V_0$ and negative definite in $V_1$;
\item[2)] the affine subspace $U=\td o+V_0$ is $\Ga$-invariant and $\Ga$ is unipotent in it;
\item[3)] the mapping of the restriction to $U$ is injective on $\Ga$ and $U/\Ga$ is strictly
causal.
\end{itemize}
Furthermore, the action of $\Ga$ in $V_1$ is linear and can be defined by an arbitrary
homomorphism $\al:\,\Ga\to\OO(\ell|_{V_1})$.
\end{theorem}

In other words, the action of $\Ga$ splits into the unipotent
affine and bounded linear ones analogously to Euclidean case. The
manifold $V/\Ga$ is isomorphic to the total space of a
topologically trivial vector bundle with the unipotent base
$U/\Ga$, the fibre $V_1$ and the holonomy representation $\al$.
The decomposition is not unique if $G$ has nontrivial fixed points
which are orthogonal to $\tau(\Ga)$. We describe two types of
$\Ga$ that classify unipotent $M=V/\Ga$ up to finite coverings.
\begin{enumerate}
\item[(\romannumeral 1)] $T,L\subset V$ are linear subspaces such
that $V=T\oplus L$, where the sum is orthogonal, $T$ is spacelike,
$\dim L=1$ and $L\cap\Int(C)\neq\emptyset$.  The group $\Ga$ is a
uniform lattice in $T$ that acts in $V$ by translations (i.e.
$\tau(\ga)=\ga$ and $\la(\ga)$ is identical in (\ref{affde})).
\item[(\romannumeral 2)] Let $v_0,v_1\in\partial C$ satisfy
$\ell(v_0,v_1)=1$ and set
\begin{eqnarray*}
&L=\bbR v_0,\quad 
W=L^\bot, \quad 
N=W\cap v_1^\bot,\quad
l_0(v)=\ell(v,v_0).
\end{eqnarray*}
The hyperplane $W=N\oplus L$ is tangent to $\partial C$ at $v_0$.
The form $\ell$ is nonpositive and degenerate in $W$, and
nondegenerate and negative in $N$.
Any $x\in N$ corresponds the following linear transformation of $V$:
\begin{eqnarray}
&\nu(x)v=v+l_0(v)x-\left(\ell(v,x)+\frac12\,l_0(v)\ell(x,x)\right)v_0.\label{nuvot}
\end{eqnarray}
A straightforward calculation shows that $\nu$ is a homomorphism of
the vector group $N$ to $\SO(\ell)$ (in fact, $\nu$ identifies it with the
factor $N$ in the Iwasawa decomposition $KAN$ for $\SO(\ell)$). Note that
for $w\in W$
\begin{eqnarray*}
&\nu(x)w=w-\ell(w,x)v_0.
\end{eqnarray*}
Further, let $T\subseteq N$ be a linear subspace and set $\td
T=T+L$. We may identify $T$ with $\td T/L$. Let $\Ga$ be a uniform
lattice in $T$ and $a$ be an $\ell$-symmetric linear mapping
\begin{eqnarray}
&a:\,T\to N,\nonumber\\
&\ell(ax,y)=\ell(x,ay), \quad x,y\in T.\label{asymm}
\end{eqnarray}
Any $a$ as above defines the affine action of $T$ in $V$ by
setting
\begin{eqnarray}
&\la(x)=\nu(ax);\label{lanua}\\
&\tau(x)=x-\frac12\,\ell(ax,x)v_0;\label{tanua}\\
&\ga_x(v)=\la(x)v+\tau(x).\nonumber
\end{eqnarray}
\end{enumerate}
Here is a necessary and sufficient condition for this action to be
free:
\begin{equation}\label{freea}
\ker({\mathbf1}+ta)=0\quad\mbox{for all}\quad t\in\bbR,
\end{equation}
where $\mathbf1$ is the identical mapping in $T$. If $a$ satisfies
(\ref{asymm}) then it admits the unique decomposition $a=a'+a''$,
where
\begin{eqnarray*}
&\quad\ell(ax,y)=\ell(a'x,y)=\ell(x,a'y), \quad x,y\in T,\\
&a':\,T\to T
\end{eqnarray*}
is the self-adjoint transformation of $T$ corresponding to the symmetric
bilinear form $\ell(ax,y)$ and $a''$ is a linear mapping
$$a'':\,T\to T^\bot\cap N.$$
The condition (\ref{freea}) can be rewritten as follows:
\begin{equation}\label{freec}
t\in\bbR,\quad a'x=tx\neq0\quad\Longrightarrow\quad a''x\neq0.
\end{equation}
Since $a'$ is self-adjoint it has real eigenvalues and
(\ref{freec}) implies that $a''\neq0$ if $a'\neq0$; then $\td
T\neq W$. Moreover, (\ref{freea}) is true if and only if the
action of $\Ga$ is free and proper (Lemma~\ref{lastl}). Also, note
that $\Ga$ consists of pure translations if $a=0$; then it is of
the type (\romannumeral 1) in fact.

\begin{theorem}\label{reduc}
If $M$ is strictly causal and unipotent then it admits a finite
covering by $V/\Ga$ with $\Ga$ as in (\romannumeral 1) or
(\romannumeral 2); in the latter case (\ref{freea}) holds.
Conversely, if $\Ga$ is as  in (\romannumeral 1) or in
(\romannumeral 2) with fulfilled (\ref{freea}) then $V/\Ga$ is
strictly causal. In both cases, $V/\Ga$ is homeomorphic to the
product of the  torus $T/\Ga$ and the vector space $V/T$.
\end{theorem}
The affine structure in $V/\Ga$ is not the product one if
$a\neq0$.

Here is  the simplest example of a non-elliptic manifold of the
type (\romannumeral2) in the least possible dimension. The group
$\Ga$ in it is a subgroup of some group of the paper \cite{Fried}.
For its generic point, either the past or the future is not closed
and contains a lightlike straight line.
\begin{example}\label{newun}\rm
Let $\dim V=4$. Set
\begin{eqnarray}
&\ell(v,v)=2v_0v_3-v^2_1-v_2^2,\qquad v=(v_0,v_1,v_2, v_3)^\top\in V=\bbR^4;\nonumber\\
&\la(n)=\left(\begin{array}{cccc}
1&0&n&\frac{n^2}2\\
0&1&0&0\\
0&0&1&n\\
0&0&0&1
\end{array}\right),\qquad\tau(n)=
\left(\begin{array}{c}
0\\
n\\
0\\
0
\end{array}\right),\qquad n\in\bbZ;\label{defex}\\
&\ga_n(v)=\la(n)v+\tau(n). \label{defon}
\end{eqnarray}
Then $\Ga=\{\ga_n\}_{n\in\bbZ}$ is the cyclic infinite group; if
$e_0,e_1,e_2,e_3$ is the standard basis of $V$ then $\Ga=\bbZ
e_1$, $l_0(v)=v_3$, $ae_1=e_2$. A calculation shows that
\begin{eqnarray}
&\ell(\ga_n(v)-v,\ga_n(v)-v)= -n^2(1+v_3^2).\label{caega}
\end{eqnarray}
Hence $\Ga$ is causal and free; clearly, $\Ga$ is proper. For
$u=(u_0,u_1,u_2,u_3)^\top$,
\begin{eqnarray}
&\ell(\ga_n(v)-u,\ga_n(v)-u)=a(v)+2nb(v)-n^2(1+u_3v_3),\label{seega}\\
&b(v)=u_2v_3-u_3v_2+u_1-v_1,\label{lofxx}\\
&a(v)=\ell(v-u,v-u).\nonumber
\end{eqnarray}
Let $\ka:\,V\to V/\Ga$ be the quotient mapping and denote
$M=V/\Ga$, $p=\ka(u)$. If $u_3=0$ then the coefficient at $n^2$ is
constant and negative. Hence
\begin{eqnarray*}
\ell(\ga_n(v)-u,\ga_n(v)-u)<0
\end{eqnarray*}
for all sufficiently large $n$. Consequently, each orbit of $\Ga$
has a finite number of common points with $C$ or $-C$. Moreover,
for any $v\in V$ there exist a neighbourhood $U_v$ of $v$ and
$n_v\in\bbN$ such that
\begin{eqnarray*}
\card(\Ga y\cap C)\leq n_v\quad\mbox{\rm for all}\quad y\in U_v.
\end{eqnarray*}
Therefore, $F_p$ and $P_p$ are closed in $M$ for all points $u$ in
the hyperplane $u_3=0$.

Let $u_3>0$. Similar arguments show that the projection of the set
$$\left\{v\in V:\,v_3>-\frac1{u_3},\  \ell(v,v)\geq0\right\}$$
to $M$ is locally closed. Hence $F_p$ is closed and $P_p\cap U_p$
is closed in $U_p$ for some neighbourhood $U_p$ of $p$ in $M$.
According to (\ref{seega}), if
\begin{eqnarray}\label{noncl}
v_3<-\frac1{u_3}
\end{eqnarray}
then $\ell(\ga_n(v)-u,\ga_n(v)-u)>0$ for all sufficiently large
$n$. Hence $\ka(v)\in P_p$  for each $v$ that satisfies
(\ref{noncl}). Clearly, $\Ga$ does not change the coordinate
$v_3$;  a computation with (\ref{lofxx}) shows what the same is
true for $b(v)$ if $1+v_3u_3=0$. Taken together with
(\ref{seega}), this implies that the set
\begin{eqnarray*}
I_u=\left\{v\in V:\,v_3=-\frac1{u_3},\ b(v)=0,\
\ell(v,v)<0\right\}
\end{eqnarray*}
is $\Ga$-invariant. Hence $\ka(I_u)$ has no common point with
$P_p$. Clearly, $I_u\neq\emptyset$ and $\ka(I_u)\subset\clos P_p$.
Thus $P_p$ is not closed. Moreover, $P_p$ contains lightlike lines
parallel to $L=\bbR e_0$.

If $u_3<0$ then, repeating this with minor changes, we get that
the past is closed but the future is not, and that the future
contains lightlike lines; also, that $M$ is strictly causal.

Let $M^0,M^+,M^-$ be the three subsets of $M$ considered above, in
the same order. Points in $M^0$ are distinguished by any of the
following properties: 1) both the past and the future are closed;
2) neither the past nor the future contain lightlike lines. Thus
$M$ is not homogeneous. The symmetry
$(v_0,v_1,v_2,v_3)\to(-v_0,v_1,-v_2,-v_3)$ commutes with $\Ga$.
Hence it defines an involution $\iota$ in $M$; evidently, $\iota$
transposes $M^+$ and $M^-$. \qed
\end{example}
In dimension 4, which is of course of particular interest,
Theorem~\ref{reduc} makes possible to find all the considered
manifolds. We omit the elliptic case which is clear. We say that
Lorentzian manifolds $(M_1,\ell_1)$ and $(M_2,\ell_2)$ are {\it
homothetic} if there exists $t>0$ such that $(M_1,t\ell_1)$ and
$(M_2,\ell_2)$ are isometric.
\begin{theorem}\label{fourc}
If a 4-dimensional strictly causal flat complete Lorentzian
manifold is not elliptic then it is homothetic to the manifold of
Example~\ref{newun}.
\end{theorem}

Theorems~\ref{split} and~~\ref{reduc}  contain the following
assertion which we formulate separately as a proposition (in fact,
it is the major point of the proof).
\begin{proposition}\label{virta}
The fundamental group $\pi_1(M)$ of a strictly causal flat
complete Lorentzian manifold $M$ is virtually abelian.
\end{proposition}

We say that $M$ is {\it elliptic} if $\Ga$ contains a finite index
subgroup whose restriction to the affine space $U$ of
Theorem~\ref{split} consists of translations; otherwise, $\Ga$ and
$M$ are said to be {\it parabolic} (the hyperbolic case is
impossible, see Lemma~\ref{hyper}). Thus elliptic manifolds
correspond to (\romannumeral 1) and parabolic to (\romannumeral 2)
with $a\neq0$ (modulo Theorem~\ref{split}). We formulate some
properties which distinguish these classes. A noncompact manifold
$M$ has {\it two ends} if for any compact set $K\subset M$ the
number of unbounded components in $M\setminus K$ is not greater
than $2$ and is equal to $2$ for some $K$. A submanifold of $M$ is
called {\it spacelike} if the restriction of $-\ell$ to it is a
Riemannian metric. A {\it Cauchy hypersurface} in a Lorentzian
manifold is a smooth spacelike submanifold of codimension 1 which
disconnects $M$ and such that each unextendible timelike (in
particular, lightlike) curve intersects it by a single point. The
existence of a Cauchy hypersurface is the main ingredient of the
definition of a globally hyperbolic Lorentzian manifold (see
\cite{HoEl}, \cite{Beem},\cite{Barbot}).
\begin{proposition}\label{disti}
Let $M$ be a strictly causal flat complete Lorentzian manifold.
Each of following conditions implies that $M$ is elliptic:
\begin{itemize}
\item[1)] $M$ is affinely equivalent to some Euclidean space form;
\item[2)] the past of any point contains no lightlike straight
line; \item[3)] the future of any point contains no lightlike
straight line;
\item[4)] $M$ admits a Cauchy hypersurface;
\item[5)] $M$ has two ends.
\end{itemize}
Conversely, elliptic manifolds satisfy 1)--4).
\end{proposition}

Causal manifolds $\bbM_4/\Ga$ with two ends, where $\Ga$ is a
discrete subgroup of a simply transitive group of automorphisms of
$\bbM_4$, were described in \cite{Fried}. Most of them are not
strictly causal. The case 4) is a direct consequence of the
paper \cite{Barbot}. We note that the past of some point in the
Cauchy hypersurface contains lightlike lines if $M$ is not
elliptic (this is a contradiction).

Let $\ga\in\Ga$, $g=\la(\ga)$. We denote by $F_\ga$ or $F_g$ the set of all fixed points
of $g$,
\[V_{\ga}=\{\la(\ga)v-v:\,v\in V\},\]
$\Int(X)$ is the interior of a set $X$, $\sspan X$ is its linear span,
$\LL(V)$ denotes the space of all linear mappings $V\to V$.
The projection of $C$ to $\bbR\bbP^{n-1}$ may be identified with the
closed unit ball $B_{n-1}$ in $\bbR^{n-1}$; the group $\OO(\ell)$ acts in it by
Moebius transformations. If $g\in\OO(\ell)$ has a fixed point in
the open unit ball then $g$ is called {\it elliptic}; this is equivalent to
\begin{equation}\label{eldef}
F_g\cap\Int(C)\neq\emptyset.
\end{equation}
If $g$ is not elliptic then
it has one or two fixed points in the unit sphere and is called {\it parabolic} or
{\it hyperbolic}, respectively.  If $g\in\OO(\ell)$ is hyperbolic and keeps $C$ then its
eigenvectors corresponding to fixed points in the sphere belong to the boundary of
$C$ and the eigenvalues are positive. If $g$ is not hyperbolic then it has no real
eigenvalues different from $\pm1$.
Let
$$F=\bigcap_{\,g\in G}F_g$$
be the set of $G$-fixed points in $V$. If
\begin{equation}\label{ellip}
F\mathbin\cap\Int(C)\neq\emptyset
\end{equation}
then we say that $G$, $\Ga$, $M$ are {\it elliptic} (for $M$, this
definition is equivalent to that was given above).

Vector spaces are always finite dimensional and real. We refer to
\cite{Wolf} for the exposition of space forms of symmetric spaces.

\section{Proof of results}
In following lemmas, $M=V/\Ga$ is assumed to be causal; the
assumption that $M$ is strictly causal is always stated
explicitly. Some of these lemmas are known and are proved for the
sake of completeness. The condition (\ref{causa}) is evidently
equivalent to
\begin{eqnarray}\label{caust}
(V_\ga+\tau(\ga))\mathbin\cap C=\emptyset\quad\mbox{for all}\ \ga\in\Ga\setminus\{e\}.
\end{eqnarray}
This implies a necessary condition for (\ref{causa}):
\begin{eqnarray}\label{causn}
V_\ga\mathbin\cap\Int(C)=\emptyset,\quad\ga\in \Ga.
\end{eqnarray}
To prove it, note that $v\in\Int(C)$ implies $tv+\tau(\ga)\in\Int(C)$
for sufficiently large $t>0$. The lemma below contains one more reformulation
of (\ref{causa}).
\begin{lemma}\label{mcaus}
$M$ is causal if and only if
\begin{equation}\label{negae}
\ell(\ga(v)-v,\ga(v)-v)<0\quad\mbox{for all}\ v\in V,\ \ga\in\Ga\setminus\{e\}.
\end{equation}
\end{lemma}
\begin{proof}
The inequality in (\ref{negae}) holds if and only if
$\ga(v)-v\notin C\mathbin\cup(-C)$. If $\ga(v)-v\in-C$ then
$\ga^{-1}(\ga(v))-\ga(v)\in C$. Hence (\ref{caust}) is equivalent
to (\ref{negae}).
\end{proof}
In the following lemma we collect some elementary facts concerning
groups of linear transformations of $C$ (they hold for any pointed
generating convex cone).
\begin{lemma}\label{eleco}
The following assertions are equivalent:
\begin{itemize}
\item[1)] $G$ is elliptic;
\item[2)] $G$ is bounded in $\GL(V)$ (i.e. its closure is compact);
\item[3)] all $G$-orbits in $\Int(C)$ are bounded;
\item[4)] there exists a bounded $G$-orbit in $\Int(C)$.
\end{itemize}
Furthermore, if $G$ is abelian then it is elliptic if and only if each
its element is elliptic.
\end{lemma}
\begin{proof}
If $v\in\Int(C)$ is $G$-fixed then the set $(C-v)\mathbin\cap(v-C)$ is a bounded
$G$-invariant neighbourhood of zero. Hence $G$ is bounded and 2) follows from 1).
Implications $2)\ \Rightarrow\ 3)\ \Rightarrow\ 4)$ are obvious.
Suppose that $Gv$ is bounded for some $v\in\Int(C)$. Then each $G$-orbit in the linear
span $W$ of $Gv$ is bounded. Hence there exists a bounded $G$-invariant neighborhood
of zero in $W$. Therefore, the restriction of $G$ to $W$ is bounded in $\GL(W)$.
Clearly, the cone $C\cap W$ is $G$-invariant. Using the averaging over the closure of
$G$ in $\GL(W)$, we get a $G$-fixed point in $\Int(C)$. Thus 4) implies 1).

If $G$ is abelian and each $g\in G$ is elliptic then the complexification of $V$
decomposes into the direct sum of $G$-eigenspaces (note that each $g$-eigenspace
is $G$-invariant and each $g\in G$ is semisimple). Since $g$-eigenvalues are bounded,
this implies that $G$ is bounded. The converse is clear.
\end{proof}
\begin{lemma}\label{hyper}
$G$ contains no hyperbolic element.
\end{lemma}
\begin{proof}
Let $\ga\in\Ga$ be hyperbolic, $g=\la(\ga)$. Then there exist $v_1,v_2$
such that
\[
gv_k=\mu_kv_k,\quad\mu_k\neq1,\quad v_k\in\partial C,\quad k=1,2,\ v_1\neq v_2.
\]
Their two-dimensional linear span $W$ intersects $\Int(C)$. Since
$W$ is $g$-invariant and $1$ is not the eigenvalue of $g$ in $W$
we have $W\subseteq V_\ga$, contradictory to (\ref{causn}).
\end{proof}

\begin{lemma}\label{cruci}
Let $V$ be a real vector space and $K$ be a subgroup of $\GL(V)$ such that
the function $\Tr $ is bounded on $K$. Then $K$ is either bounded or
reducible.
\end{lemma}
\begin{proof}
Set
\begin{eqnarray*}
&A=\sspan K,\\
&B=\{x\in\LL(V):\,\sup_{h\in K}|\Tr xh|<\infty\},\\
&N=A^\bot=\{x\in\LL(V):\,\Tr xh=0\quad\mbox{for all}\quad h\in K\}.
\end{eqnarray*}
Due to the assumption of the lemma, $K\subset B$. Hence
$A\subseteq B$. Clearly, $N\subseteq B$ and $A$ is an algebra. Set
$I=A\mathbin\cap N$; standard arguments shows that $I$ is a
two-side ideal in $A$.

If $I=0$ then the bilinear form $\Tr xy$ is nondegenerate in $A$.
Therefore, any linear functional on $A$ is bounded on $K$. Thus
$K$ is bounded.

For any $x\in I$ and positive integer $n$, $\Tr x^n=0$; hence $x$
is nilpotent, and Engel's theorem implies that the space $Z=\{v\in
V:\,Iv=0\}$ is nontrivial. If $I\neq0$ then $Z\neq V$. Since $Z$
is evidently $A$-invariant, $K$ is reducible.
\end{proof}
\begin{lemma}\label{fixed}
Let $K$ be a subgroup of $\OO(\ell)$ that keeps $C$ and contains no
hyperbolic element. Then $K$ has a fixed point in $C\setminus\{0\}$.
\end{lemma}
\begin{proof}
The assumption implies that eigenvalues of any $h\in K$ are
contained in the unit circle. Hence $|\Tr h|\leq\dim V$. If $K$ is
irreducible then it is bounded by Lemma~\ref{cruci}. Then its
closure is compact and we get $K$-fixed points in $\Int(C)$ by
averaging.

Let $K$ be reducible. Then there exists a proper $K$-invariant
space $W$. The invariant space $L=W\mathbin\cap W^\bot$ is at most
one dimensional since the Lorentzian form $\ell$ vanishes on it.
If $\dim L=1$ then all points of the line $L$ are $K$-fixed:
otherwise, $K$ either contains an element that has positive
eigenvalue $\mu\neq1$ which is hyperbolic or does not keep $C$
(note that one of the two halflines in $L$ lie in $\partial C$).
Thus the assertion is true if $L\neq0$. If $L=0$ then $V=W\oplus
W^\bot$. Since $\ell$ is nondegenerate in each of these spaces,
either one of them is one dimensional and intersects $\Int(C)$ or
the restriction of $\ell$ to one of them is Lorentzian. The first
case is clear. Thus we may use the induction on $\dim V$ starting
with the obvious case $\dim V=1$.
\end{proof}
\begin{lemma}\label{eqria}
If $G$ is elliptic then $M$ admits a flat Riemannian metric making it an Euclidean space
form (i.e. a flat geodesically complete Riemannian manifold) such that the identical mapping
is affine.
\end{lemma}
\begin{proof}
By definition, there exists $v_0\in F\mathbin\cap\Int(C)$. The form
\begin{equation}\label{semid}
t\ell(v_0,v)^2-\ell(v,v),
\end{equation}
is $G$-invariant. Moreover, for sufficiently large $t>0$ it is positive definite. Hence
it induces the desired metric on $M=V/\Ga$.
\end{proof}
In that follows, we assume that $G$ is not elliptic if the
contrary is not stated explicitly. By Lemma~\ref{fixed}, there
exists
\begin{eqnarray}\label{fixv0}
v_0\in\partial C\cap F.
\end{eqnarray}
For this vector $v_0$, we define $L$, $W$, $l_0$ as in
({\romannumeral2}). By (\ref{eldef}), $F\cap\Int(C)=\emptyset$.
Hence
\begin{eqnarray}
L\subseteq F\subseteq W.
\end{eqnarray}
Since $l_0$ is $G$-invariant,
\begin{eqnarray*}
l_0(\ga(v)-v)=l_0(\tau(\ga));
\end{eqnarray*}
in particular, the left part does not depend on $v$. The following lemma is obvious.
\begin{lemma}\label{fahyp}
The group $\Ga$ keeps the family of hyperplanes in $V$ that are
parallel to $W$; moreover, $G$ and the commutator group of $\Ga$ keep each of them.
The mapping
\begin{equation}\label{gatau}
\al:\,\ga\to\l_0(\tau(\ga)).
\end{equation}
is a homomorphism $\Ga\to\bbR$.\qed
\end{lemma}
The aim of subsequent lemmas is to prove that
\begin{eqnarray}\label{gainh}
\tau(\Ga)\subset W.
\end{eqnarray}
\begin{lemma}\label{paral}
If $V_\ga\supseteq L$ then $\tau(\ga)\in W$.
\end{lemma}
\begin{proof}
Let $v\in V$ be such that $\la(\ga)v-v=v_0$. Set $b=\tau(\ga)$. For all $t\in\bbR$,
\begin{eqnarray}
\ell(\ga(tv)-tv,\ga(tv)-tv)=\ell(tv_0+b,tv_0+b)=2tl_0(b)+\ell(b,b).
\end{eqnarray}
If $b\notin W$ then $l_0(b)\neq0$. Hence $\ell(\ga(tv)-tv)>0$ for some $t\in\bbR$
and $M$ cannot be causal by Lemma~\ref{mcaus}. Thus $b\in W$.
\end{proof}
Note that $V_\ga\perp v_0$ since $v_0$ is $G$-fixed. Hence
$V_\ga\subseteq W$ and
\begin{eqnarray}\label{equcl}
V_\ga\mathbin\cap C=\{0\}\quad\Longleftrightarrow\quad V_\ga\mathbin\cap L=\{0\}.
\end{eqnarray}
\begin{lemma}\label{ellip}
$V_\ga\mathbin\cap L=\{0\}$ if and only if $\ga$ is elliptic.
\end{lemma}
\begin{proof}
Since $V_\ga\subseteq W$ and $W^\bot=L$,  the assumption $V_\ga\mathbin\cap L=\{0\}$
implies that $\ell$ is negative definite on $V_\ga$. Therefore,
\begin{eqnarray*}
V=V_\ga\oplus V_\ga^\bot
\end{eqnarray*}
and $\la(\ga)$ generates a bounded subgroup of $\GL(V_\ga)$.
Clearly, $\la(\ga)$ is identical in $V_\ga^\bot$. Thus $\la(\ga)$
generates a bounded subgroup of $\GL(V)$ and is elliptic by
Lemma~\ref{eleco}. Conversely, if $\ga$ is elliptic then
$\la(\ga)$ is semisimple according to the same Lemma. Since
$L\subseteq F_\ga$, this implies $V_\ga\mathbin\cap L=\{0\}$.
\end{proof}
\begin{corollary}\label{parak}
$V_\ga\supseteq L$ if and only if $\ga$ is parabolic.
\end{corollary}
\begin{proof}
Since $\dim L=1$, $V_\ga\not\supseteq L$ is equivalent to $V_\ga\mathbin\cap L=\{0\}$.
\end{proof}
\begin{corollary}\label{ellga}
If $\tau(\ga)\notin W$ then $\ga$ is elliptic.
\end{corollary}
\begin{proof}
Combine Lemma~\ref{paral} and Corollary~\ref{parak}.
\end{proof}
Set
\begin{eqnarray}\label{gapri}
\Ga'=\{\ga\in\Ga:\,\tau(\ga)\in W\},\quad G'=\la(\Ga').
\end{eqnarray}
By Lemma~\ref{fahyp}, $\Ga'$ and $G'$ are normal in $\Ga$, $G$, respectively.
\begin{lemma}\label{gprig}
If $G'$ is elliptic then $G$ is elliptic.
\end{lemma}
\begin{proof}
Let $F'$ be the set of all $G'$-fixed points in $V$. If $G'$ is elliptic then
\begin{eqnarray}\label{gprim}
F'\mathbin\cap \Int(C)\neq\emptyset.
\end{eqnarray}
Since $G'$ is normal, $F'$ is $G$-invariant. The action of $G$ in $F'$
can be considered as the action of $G/G'$. By Lemma~\ref{fahyp}, $G/G'$ is
abelian. Due to Corollary~\ref{ellga}, each $g\notin G'$ is elliptic.
Taken together with (\ref{gprim}), this satisfies the assumption of
Lemma~\ref{eleco} for the group $G/G'$, the space $F'$, and the cone $F'\cap C$.
Therefore, $G$ has a fixed point in $\Int(C)\mathbin\cap F'\subseteq\Int(C)$.
Hence $G$ is elliptic.
\end{proof}
According  to Lemma~\ref{eqria} and Lemma~\ref{gprig}, there are
two levels of the problem: first, the case of non-elliptic $G$
under the additional assumption $\Ga=\Ga'$, and second, the
covering $\mu:\,M'\to M$, where $M'=V/\Ga'$. We show that the
second is trivial (i.e. $\Ga=\Ga'$).

Condition (\ref{gainh}) and $G$-invariance of $L$ define the
affine action of $\Ga$ in the quotient space $W/L$ which may be
identified with $N$. Since $W=L^\bot$, the form $\ell$ induces a
negative definite $G$-invariant form on $W/L$ which coincides with
its restriction to $N$. Let $\phi:\,V\to V/L$ and $\ka:\,V\to
M=V/\Ga$ be the quotient mappings. Then $\phi(W)=N$.
\begin{lemma}\label{parab}
Let (\ref{gainh}) be true. If $M$ is causal then the action of
$\Ga$ in $N=W/L$ is free; if $M$ is strictly causal then it is
discontinuous (i.e. each orbit is discrete).
\end{lemma}
\begin{proof}
It follows from (\ref{causa}) that the projection $\phi:\,W\to N$
is one-to-one on each $\Ga$-orbit. This proves the first
assertion.

Let $u,w\in W$ and $\{\ga_n\}$ be a sequence in $\Ga$ such that
\begin{eqnarray*}
\lim_{n\to\infty}\phi(\ga_n(w))=\phi(u).
\end{eqnarray*}
Then $\ga_n(w)=u+w_n+t_nv_0$, where $w_n\to0$ in $W$ as
$n\to\infty$ and $t_n\in\bbR$. Since $\ga_n(w+tv_0)=\ga_n(w)+tv_0$
for all $t\in\bbR$, this implies that $\ka(w+tv_0)$ lies in the
closure of $\ka(u+L)$ in $M$ for each $t\in\bbR$. Also, it follows
that for any neighbourhood $U$ of $p=\ka(u)$ in $M$ and some
$s\in\bbR$ the inclusion $\ka(w+sv_0)\in U$ holds. Let $U$ be such
that $U\cap F_p$ and $U\cap P_p$ are closed in $U$. Set
$L^+=\{tv_0:\,t\geq0\}$; clearly, $\ka(w+sv_0)$ belongs to the
closure of at least one of the sets $\ka(u+L^+)$ and $\ka(u-L^+)$.
On the other hand, $\ka(u-L^+)\cap U$ and $\ka(u+L^+)\cap U$ are
closed in $U$. Indeed, due to (\ref{gainh}), $\ka(W)$ is closed in
$M$ and
\begin{eqnarray*}
\ka(u+L^+)=F_p\cap\ka(W),\qquad \ka(u-L^+)=P_p\cap\ka(W)
\end{eqnarray*}
(recall that $F_p=\ka(u+C)$, $P_p=\ka(u-C)$). Therefore,
$\ka(w+sv_0)\in\ka(u+L)$; this means that $\phi(\ga(w))=\phi(u)$
for some $\ga\in\Ga$.
Hence each orbit of $\Ga$ in $W/L$ is closed.
Consequently, all orbits are discrete.
\end{proof}
\begin{corollary}\label{iseuc}
If (\ref{gainh}) is true and $M$ is strictly causal then $N/\Ga$
is an Euclidean space form; in particular, $\Ga$ is virtually
abelian.\qed
\end{corollary}

\begin{lemma}\label{fiind}
If $\Ga$ is not elliptic then any abelian subgroup
$\wtd\Ga\subseteq\Ga$ of finite index contains a parabolic
element.
\end{lemma}
\begin{proof}
If $\wtd\Ga$ consists of elliptic elements then $\wtd G=\la(\wtd\Ga)$
is elliptic by Lemma~\ref{eleco}.  Then $\wtd G$ has a fixed point $w_0\in\Int(C)$.
The orbit $Gw_0$ is finite; the sum of vectors in it belongs to $\Int(C)$ and
is $G$-fixed.
\end{proof}
Any affine action in $\bbR^m$ has the natural extension to linear one in $\bbR^{m+1}$.
Precisely, the mapping $x\to ax+b$ can be realized as the restriction to the hyperplane
$x_{m+1}=1$ of the linear mapping $\td a$ that is equal to $a$ on $\bbR^m$
(embedded to $\bbR^{m+1}$ as the hyperplane $x_{m+1}=0$) and satisfies
\begin{eqnarray}\label{liaff}
\td ae_{m+1}=e_{m+1}+b,\quad e_{m+1}=(0,\dots,0,1).
\end{eqnarray}
Thus we may consider $\Ga$ as a linear group in $\td V=V\oplus\bbR$.
\begin{lemma}\label{aflid}
Let $\Ga$ be an abelian group of affine transformations in $V$. Then there exist
$\td o\in V$ and linear subspaces $V_0,V_1\subseteq V$ such that
\begin{eqnarray}
&V=V_0\oplus V_1,\label{vdeco}
\end{eqnarray}
the affine subspace $U=\td o+V_0$ is $\Ga$-invariant, $\Ga$ is unipotent in $U$,
and $\Ga$ is linear in $V_1$.
\end{lemma}
\begin{proof}
Let $\td V=V\oplus\,\bbR$ be as above and $\pi$ be the projection to the first
component along the second one.
Since $\Ga$ is abelian, it admits the triangular realization in some linear base
of $\td V$. Diagonal elements are characters of $\Ga$ and $\td V$ decomposes over them.
Let $\td V_0$ be the component of the trivial character. Then
\begin{eqnarray}\label{jorvv}
\td V=\td V_0\oplus V_1,
\end{eqnarray}
where $V_1$ is the sum of all other components. Clearly, $V$ is $\Ga$-invariant
(we identify $V$ with $V\oplus\,0$). According to (\ref{liaff}),
the representation of $\Ga$ in $\td V/V$ is trivial. Hence
\begin{eqnarray}
&V_1\subseteq V\label{vodin}
\end{eqnarray}
and $(V\oplus1)\cap\td V_0\neq\emptyset$.
The set $U=\pi((V\oplus1)\cap\td V_0)$ is an affine subspace of $V$.
Since $\td V_0$ and $V\oplus1$ are $\Ga$-invariant, $U$ is invariant
with respect to the affine action of $\Ga$ in $V$.
Pick any $\td o\in U$ and set $V_0=\td V_0\cap V$. Then $\td o\oplus\,1\in\td V$,
$\td V_0=\td o\oplus\,1+V_0$, $U=\td o+V_0$
and (\ref{vdeco}) follows from (\ref{jorvv}) and (\ref{vodin}).
\end{proof}
We need a lemma that reduces the virtually abelian case to the abelian one.
\begin{lemma}\label{vabab}
Let a group $\Ga\subset\Aff(V)$ be  virtually abelian. Then there exists a finite index
abelian subgroup $\td\Ga\subseteq\Ga$ such that
\[
\al\in\Aff(V),\ \al^{-1}\Ga\al=\Ga\quad\Longrightarrow\quad\al^{-1}\td\Ga\al=\td\Ga.
\]
\end{lemma}
\begin{proof}
Let $\Ga'$ be an abelian subgroup of finite index in $\Ga$, and
let $\wh\Ga'$ be its algebraic closure. Then $\wh\Ga=\Ga\wh\Ga'$
is closed being a finite union of closed sets; hence $\wh\Ga$ is
the algebraic closure of $\Ga$ (in particular, $\wh\Ga$ is a
group). Let $\wh\Ga_0$ be the identity component of $\wh\Ga$. The
group $\td\Ga=\Ga\cap\wh\Ga_0$ evidently satisfies the lemma.
\end{proof}
Any affine transformation $\ga$ is the composition of its linear part
$\la(\ga)$ and the translation part $v\to v+\tau(\ga)$.
They commute if and only if $\tau(\ga)$ is $\la(\ga)$-fixed:
\begin{equation}\label{comlt}
\la(\ga)\tau(\ga)=\tau(\ga).
\end{equation}
Of course, the decomposition and (\ref{comlt}) depend on the
choice of origin: if it is removed to $v\in V$ then $\tau(\ga)$ is
replaced by $\ga(v)-v=\la(\ga)v-v+\tau(\ga).$ For the
linearization of $\ga$ in $V\oplus\bbR$,  the existence of the
origin satisfying (\ref{comlt}) is equivalent to
$\ker({\mathbf1}-\ga)\not\subseteq V\oplus0$. Groups $\Ga$ of
({\romannumeral2}) do not have this property.
\begin{lemma}\label{trudn}
Suppose that (\ref{gainh}) is not true,  $\Ga'$ is an abelian subgroup of
$\Ga$, let $V_0$, $V_1$, $\td o$ be as in Lemma~\ref{aflid} for $\Ga'$ and remove
origin to $\td o$. Then (\ref{comlt}) holds for any parabolic $\ga\in\Ga'$.
\end{lemma}
\begin{proof}
Let $\ga\in\Ga'$ be parabolic and set
$W_0=W\mathbin\cap V_0$, $g=\la(\ga)$, $b=\tau(\ga)$.
Note that
$$b\in W_0$$
by Corollary~\ref{ellga} and Lemma~\ref{aflid}.
Since $g$ is orthogonal with respect
to the form $\ell$ in $N=W/L$, it is semisimple in $N$. Hence $g$ is identical
in $W_0/L$. Thus there exists a linear functional $l$ on $W_0$ such that
\begin{equation}\label{ldefi}
gu=u+l(u)v_0
\end{equation}
for all $u\in W_0$. Clearly, $l(v_0)=0$.
If $l(b)=0$ then $gb=b$ and (\ref{comlt}) is true.
Therefore, it is sufficient  to prove that
\begin{equation}\label{nonca}
\ell(\ga^n\circ\nu(o),\ga^n\circ\nu(o))>0
\end{equation}
for some $\nu\in\Ga$ and $n\in\bbZ$ assuming that (\ref{gainh}) is false and
\begin{eqnarray}\label{nospl}
l(b)\neq0.
\end{eqnarray}
Also, we may assume that $\Ga'$ is unipotent since $V_0$ is $\Ga'$-invariant.
Let $w\in W_0$. A calculation shows that
\begin{eqnarray}\label{itera}
\ga^n(w)=w+nb+\left(nl(w)+\frac{n(n-1)}{2}\,l(b)\right)v_0.
\end{eqnarray}
If $v_1\in V_0\setminus W_0$ then $\ell(v_0,v_1)\neq0$ since $\ell$
is nondegenerate in $V_0$. Replacing $v_1$ by $t_1v_1+t_0v_0$ for suitable
$t_1,t_0\in\bbR$, we may assume
\begin{eqnarray}\label{norml}
\ell(v_1,v_1)=0,\quad
\ell(v_1,v_0)=1.
\end{eqnarray}
If $gv_1=v_1$ then
$v_0+v_1$ is a $g$-fixed point in $\Int(C)$ but $g$ is supposed to
be parabolic. Hence
\begin{equation}\label{uonez}
u_0=gv_1-v_1\neq0.
\end{equation}
Clearly, $u_0\in W$; we claim that
\begin{equation}\label{uonol}
u_0\in W\setminus L.
\end{equation}
If $u_0\in L$ then $u_0=gu_0=g^{-1}u_0$. Since $v_1$ is lightlike,
\begin{eqnarray*}
\ell(gv_1,v_1)=\ell(u_0,v_1)=\ell(g^{-1}u_0,v_1)=-\ell(g^{-1}v_1,v_1)=
-\ell(v_1,gv_1).
\end{eqnarray*}
Consequently, $\ell(u_0,v_1)=0$ if $u_0\in L$; by (\ref{norml}), $u_0=0$
but this contradicts to (\ref{uonez}). This proves (\ref{uonol}).
Since $u_0\in W$, $gu_0-u_0\in L$.
Hence $gu_0-u_0=u_0-g^{-1}u_0$; using (\ref{norml}) we get
\begin{eqnarray*}
\ell(gu_0-u_0,v_1)=\ell(g^2v_1-2gv_1+v_1,v_1)=\ell(gv_1-2v_1+g^{-1}v_1,v_1)=\\
\ell(gv_1,v_1)+\ell(g^{-1}v_1,v_1)= 2\ell(gv_1,v_1)=
-\ell(gv_1-v_1,gv_1-v_1)=-\ell(u_0,u_0).
\end{eqnarray*}
Due to (\ref{uonol}), this implies $\ell(gu_0-u_0,v_1)\neq0$.
Therefore, $gu_0-u_0\neq0$ and
\begin{eqnarray}\label{iterg}
g^nv_1=v_1+nu_0+\frac{n(n-1)}{2}pv_0,\qquad pv_0=gu_0-u_0.
\end{eqnarray}
According to (\ref{norml}), (\ref{uonol})
and the calculation above,
\begin{eqnarray}\label{pelll}
p=\ell(gu_0-u_0,v_1)=-\ell(u_0,u_0)\neq0.
\end{eqnarray}
Let $v\in V_0$, $v=tv_1+w$, where
$$t=\ell(v,v_0)$$
and $w\in W_0$. By  (\ref{itera}) and (\ref{iterg}),
\begin{eqnarray*}
\ga^n(v)=\ga^n(tv_1+w)=tg^nv_1+\ga^n(w)=tv_1+w+n(tu_0+b+l(w)v_0)+\\
\frac{n(n-1)}{2}(tp+l(b))v_0.
\end{eqnarray*}
Relations
$$\ell(u_0,v_0)=\ell(b,v_0)=\ell(v_1,v_1)=\ell(v_0,v_0)=\ell(w,v_0)=0,\quad\ell(v_1,v_0)=1,$$
imply
\begin{eqnarray}
\ell(\ga^n(v),\ga^n(v))=n^2(t(tp+l(b))+\ell(tu_0+b,tu_0+b))+O(n)=\phantom{xx}\nonumber\\
n^2(l(b)+\ell(b,b)+2t\ell(u_0,b))+O(n),\label{asiga}
\end{eqnarray}
where the latter equality holds due to (\ref{pelll}).
By (\ref{ldefi}) and (\ref{norml}), for all $u\in W_0$
\begin{eqnarray*}
l(u)=\ell(l(u)v_0,v_1)=\ell(gu-u,v_1)=\ell(u,g^{-1}v_1-v_1)=-\ell(u,g^{-1}u_0).
\end{eqnarray*}
Since $gu_0-u_0=sv_0$ for some $s\in\bbR$ and $v_0$ is $g$-fixed, we have
$g^{-1}u_0=u_0-sv_0$. Therefore,
\begin{eqnarray}\label{luuno}
l(u)=-\ell(u,u_0).
\end{eqnarray}
It follows from (\ref{nospl}) that the coefficient at $n^2$ is
positive for sufficiently large $t$ of the same sign as $-l(b)$.
Lemma~\ref{fahyp} and the assumtion $\tau(\Ga)\not\subset W$ imply
the existence of $\nu\in\Ga$ such that this property is true for
$t=\ell(\tau(\nu),v_0)$. Then (\ref{nonca}) holds for all
sufficiently large $n\in\bbZ$.
\end{proof}
The linear part $g$ of $\ga$ keeps paraboloids $\partial
C\mathbin\cap(W+x)$ whose axes are collinear to $v_0$. The
calculation above shows that the asymptotic behavior of powers
$\ga^n$ is determined by $\la(\ga)$.
\begin{lemma}\label{parfi}
Let $\Ga$ be abelian and unipotent. Suppose that (\ref{comlt})
holds for each $\ga\in\Ga$. Then
\begin{eqnarray}\label{fixgg}
\la(\ga)\tau(\nu)=\tau(\nu)
\end{eqnarray}
for all $\ga,\nu\in\Ga$. Moreover, the mapping $\ga\to\tau(\ga)$
is an isomorphism of $\Ga$ onto a uniform lattice in the linear
span $T$ of $\tau(\Ga)$.
\end{lemma}
\begin{proof}
Applying (\ref{comlt}) to the composition $\ga\circ\nu=\nu\circ\ga$, we get
\begin{eqnarray*}
\la(\ga)\la(\nu)(\la(\ga)\tau(\nu)+\tau(\ga))=
\la(\nu)\tau(\ga)+\tau(\nu).
\end{eqnarray*}
Since $\la(\ga)$ and $\la(\nu)$ commute and (\ref{comlt}) holds
for $\ga$ and $\nu$, the left side is equal to
$\la(\ga)^2\tau(\nu)+\la(\nu)\tau(\ga)$. Therefore,
$\la(\ga)^2\tau(\nu)=\tau(\nu)$. Since $\la(\ga)$ is unipotent,
$\la(\ga)$ and $\la(\ga)^2$ have identical sets of fixed points.
This proves (\ref{fixgg}). The composition law for affine mappings
and (\ref{fixgg}) imply
$$\tau(\ga\circ\nu)=\tau(\ga)+\tau(\nu).$$
Since $\Ga$ acts discontinuously and freely, $\tau(\Ga)$ is a
uniform lattice in $T$ and the mapping $\tau$ is an isomorphism.
\end{proof}
For $T$ and $\Ga$ as in Lemma~\ref{parfi}, we may identify $\Ga$
and $\tau(\Ga)$. In other words, we may assume that $\Ga$ is
embedded to $T$ as a uniform lattice and acts in $V$ by
$x\to\ga_x$, where
\begin{eqnarray}\label{tauid}
\ga_x(v)=\la(x)v+x
\end{eqnarray}
for any $x\in\Ga$.
\begin{lemma}\label{norga}
Let $\Ga$, $T$  be as above.
If the affine mapping $\nu(v)=cv+b$ normalizes $\Ga$ then either $b\in W$ or
$\la(x)b=b$ for all $x\in\Ga$.
\end{lemma}
\begin{proof}
By the assumption, for all $x\in\Ga$
$$\nu^{-1}\circ\ga_x\circ\nu=\ga_{\si(x)},$$
where $\si$ is an automorphism of $\Ga$. A computation with
(\ref{tauid}) shows that
\[
\si(x)=\tau(\nu^{-1}\circ\ga_x\circ\nu)=-c^{-1}((\la(x)b-b)+x).
\]
If $b\notin W$ and $\la(x)b\neq b$ for some $x\in\Ga$ then calculations
of Lemma~\ref{trudn} (precisely, (\ref{iterg}) and (\ref{pelll})) show
that $\la(nx)$ is quadratic on $n$.  Then $\si(nx)$ is also quadratic on $n$
but each automorphism of $\Ga$ extends to a linear transformation of $T$.
\end{proof}
\begin{proposition}\label{betow}
If $M=V/\Ga$ is strictly causal and non-elliptic then $\tau(\Ga)\subset W$.
\end{proposition}
\begin{proof}
Suppose that there exists $\nu\in\Ga$, $\nu(v)=cv+b$, such that
$b\notin W$. Let $\Ga'$ be defined by (\ref{gapri}). Clearly,
$V/\Ga'$ is strictly causal and $\Ga'$ is normal in $\Ga$. By
Corollary~\ref{iseuc}, $\Ga'$ is virtually abelian. Let $\Ga''$ be
the abelian subgroup of finite index in $\Ga'$ with the property
of Lemma~\ref{vabab}. Then $\Ga''$ is normal in $\Ga$ and has
finite index in $\Ga'$. Due to Lemma~\ref{aflid}, removing the
origin if necessary, we may assume that $V=V_0\oplus V_1$, where
$V_0$ and $V_1$ are linear subspaces, the decomposition is
$\Ga''$-invariant, $\Ga''$ acts linearly in $V_1$ and affinely in
$V_0$. Standard arguments show that $V_0$ is
$\Ga$-invariant. 
Hence $b\in V_0$ and we may assume that $\Ga''$ is unipotent.
Lemma~\ref{trudn} implies (\ref{comlt}) for any parabolic
$\ga\in\Ga''$. If $\la$ is elliptic and unipotent then
$\la(\ga)={\mathbf1}$. Hence (\ref{comlt}) is true for all
$\ga\in\Ga''$. Applying consequently Lemma~\ref{parfi} and
Lemma~\ref{norga}, we get
$$\la(\ga)b=b$$
for all $\ga\in\Ga''$. Since $b\notin W$, the two dimensional
plane passing through $b$ and $v_0$ intersects $\Int(C)$. Thus
$\ga$ is elliptic. Hence $\Ga''$ consists of elliptic elements; by
Lemma~\ref{eleco}, $\Ga''$ is elliptic. Then $\Ga'$ and $\Ga$ are
elliptic by Lemma~\ref{fiind} and Lemma~\ref{gprig}, respectively.
\end{proof}
\begin{proof}[Proof of Proposition~\ref{virta}]
If $\Ga$ is elliptic then the assertion holds due to
Lemma~\ref{eqria}. If $\Ga$ is not elliptic then (\ref{gainh}) is
true for $\Ga$ by Proposition~\ref{betow}. Hence $\Ga$ is
virtually abelian by Corollary~\ref{iseuc}.
\end{proof}
\begin{proof}[Proof of Theorem~\ref{split}]
Due to Proposition~\ref{virta}, we may assume that $\Ga$ is
abelian. Then Lemma~\ref{aflid} implies the existence of the
decomposition; moreover, 2) is an easy consequence of the
construction (recall that $V_1$ is the sum of all nontrivial
components and $V_0$ is the unipotent component of $G$ in $V$).
Further, 3) holds since the action of $\Ga$ is free.  The
decomposition is evidently $\ell$-orthogonal and $G$-invariant. By
Lemma~\ref{fixed}, $C$ contains a nontrivial fixed point $v_0$.
Since $\dim v_0^\bot\cap C\leq1$, $V_1$ is spacelike. Hence $\ell$
is negative definite in it. This implies that $\ell$ is Lorentzian
in $V_0=V_1^\bot$. Remaining assertions are clear.
\end{proof}
\begin{lemma}\label{reafo}
Let $\Ga$ be unipotent.  Suppose that  $\tau(\Ga)\subset W$ and
$\phi\circ\tau$ isomorphically embeds $\Ga$ to the vector group
$W/L$ as a uniform lattice in some its linear subspace. Then there
exists $v_1\in\partial C$ and subspaces $T\subseteq N\subseteq W$
such that the action of $\Ga$ in $V$ is subject to formulas of
(\romannumeral2).
\end{lemma}
\begin{proof}
Note that any choice of $v_1\in\partial C$ (we shall refine it
later) such that $\ell(v_0,v_1)=1$ identifies $W/L$ with the space
$N=W\cap v_1^\bot$, the restriction of $\phi$ to $W$ with the
orthogonal projection $\pi$ onto $N$, and $T$ with the linear span
of $\phi(\tau(\Ga))$. We may assume that $\Ga$ is embedded to $T$
by $\phi\circ\tau$; then $\phi\circ\tau(x)=x$ for all $x\in\Ga$.
Since $\la(x)$ is $\ell$-orthogonal and unipotent in $W/L$ for any
$x\in\Ga$, $\la(x)=\mathbf1$ in $W/L$. Thus
\begin{eqnarray}
&\tau(x)=x+\xi(x)v_0,\nonumber\\
&\la(x)w=w+\eta(x,w)v_0,\label{laide}
\end{eqnarray}
where $\xi$ and $\eta$ are real valued functions, $x\in\Ga$ and
$w\in W$.
Clearly, $\eta(x,w)$ is linear on $w\in W$. Since $v_0$ is
$G$-fixed and
\begin{eqnarray*}
&\tau(x+y)=\la(x)\tau(y)+\tau(x)=\la(y)\tau(x)+\tau(y),
\end{eqnarray*}
we get
\begin{eqnarray*}
&\la(x)\tau(y)=\la(x)(y+\xi(y)v_0)=y+(\eta(x,y)+\xi(y))v_0,\\
&\tau(x+y)-x-y=(\eta(x,y)+\xi(x)+\xi(y))v_0=\xi(x+y)v_0.
\end{eqnarray*}
This implies $\eta(x,y)=\xi(x+y)-\xi(x)-\xi(y)$ for $x,y\in T$.
Hence $\eta$ is symmetric bilinear on $T$ and
$$\xi(x)=\frac12\,\eta(x,x)+\ze(x),$$
where $\ze$ is a linear functional on $T$. Let $z\in T$ be such
that $\ze(x)=\ell(z,x)$ for all $x\in T$. Replacing $T$ by the
space $\{x+\ze(x)v_0:\,x\in T\}$, $v_1$ by a suitable combination
of $v_1-z$ and $v_0$, and redefining $N$, we may assume that
$\ze=0$. For $x\in \Ga$, set
$$ax=\pi(\la(x)v_1),$$
where $\pi$ is the orthogonal projection onto $N$. Then
$$
\la(x)v_1=\si(x)v_1+ax+\om(x)v_0,
$$
where $\si,\om$ are functions on $\Ga$. Since $ax\in N$,
$v_1,v_0\in\partial C$, $\ell(v_1,v_0)=1$, $\la(x)\in\SO(\ell)$,
and $v_0$ is $G$-fixed,
\begin{eqnarray*}
&\si(x)=\ell(\la(x)v_1,v_0)=\ell(v_1,\la(x)^{-1}v_0)=1;\\
&0=\ell(\la(x)v_1,\la(x)v_1)=2\om(x)+\ell(ax,ax).
\end{eqnarray*}
If $w\in W$ then $w-\la(x)w\in L$. Hence $\pi(\la(x)w)=\pi(w)$; in
particular, this is true for $w=v_1-\la(y)v_1$ and  we get
\begin{eqnarray*}
&a(x+y)=\pi(\la(x+y)v_1)=\pi(\la(x)v_1)+\pi(\la(x)(\la(y)v_1-v_1))
=ax+ay,\\
&ao=\pi(v_1)=o
\end{eqnarray*}
for all $x,y\in\Ga$. Therefore, $a$ extends to a linear operator
$T\to N$. Clearly, $\la(x)^{-1}=\la(-x)$. Due to (\ref{norml}) and
(\ref{laide}), if $x\in\Ga$ and $w\in N$ then
\begin{eqnarray*}
\eta(x,w)=\ell(v_1,\la(x)w)=\ell(\la(-x)v_1,w)=-\ell(ax,w).
\end{eqnarray*}
Since $\la(x)v_0=v_0$, (\ref{laide}) extends the derived formula
to the case $w\in W$ by setting $\eta(x,v_0)=0$. For any $v\in V$
we have $v=l_0(v)v_1+w$, where $w\in W$. Combining equalities
above, we get formulas of ({\romannumeral2}) for $\la$ and $\tau$:
\begin{eqnarray}
&\la(x)v
=v+l_0(v)ax-(\ell(ax,v)+\frac12\,l_0(v)\ell(ax,ax))v_0,\label{forla}\\
&\tau(x)=x-\frac12\,\ell(ax,x)v_0,\label{forta}
\end{eqnarray}
where $v\in V$, $x\in\Ga$. Since $\eta$ is symmetric, $a$ is
$\ell$-symmetric in $T$.
\end{proof}

\begin{lemma}\label{lastl}
Each of following conditions is equivalent to (\ref{freea}):
\begin{itemize}
\item[1)] the action ({\romannumeral2}) of $T$ in $V$ is free;
\item[2)] the action ({\romannumeral2}) of $\Ga$ in $V$ is free
and proper.
\end{itemize}
\end{lemma}
\begin{proof}
The action of $T$ in $V$ is not free if and only if
\begin{eqnarray*}
\ga_x(v)-v=(x+l_0(v)ax)-
\left(\ell(v,ax)+\frac12\,l_0(v)\ell(ax,ax)+\frac12\,\ell(ax,x)\right)v_0=0
\end{eqnarray*}
for some $x\in T\setminus\{0\}$, $v\in V$. This equality is
equivalent to the system
\begin{eqnarray*}
\left\{
\begin{array}{r}
x+l_0(v)ax=0;\cr \ell(v,ax)=0.
\end{array}
\right.
\end{eqnarray*}
The first equation has a solutions $x,v$ such that $x\neq0$ if and
only if (\ref{freea}) is not true.  If $\dim V>1$ then vector $v$
can be removed without changing $l_0(v)$ to satisfy the second
equation (note that $ax\notin L\setminus\{0\}$). Since the
assertion is clear for $\dim V=1$, (\ref{freea}) is the same as
1).

Hyperplanes $X_s=\{v\in V:\,l_0(v)=s\}$ are $T$-invariant; the
space $L$ is parallel to each of them. Hence $T$ acts in
hyperplanes $\phi(X_s)$ in $V/L$ by translations $v\to v+x+sax$.
Therefore, the action of the uniform lattice $\Ga\subset T$ in all
$\phi(X_s)$ is free and proper if and only if $T$ acts freely. If
$\Ga$ is free and proper in all these hyperplanes then it is free
and proper in $V/L$ since the fundamental parallelepiped for
lattices in $\phi(X_s)$ depends on $s$ continuously. Then the same
is true for $V$. If $\Ga$ is not free in $V/L$ then it is not free
in $V$ due to the calculation above. A minor modification of these
arguments shows that the action of $\Ga$ is not proper if $T$ is
not free in $V/L$ but $\Ga$ is free (then there exists $v\in V$
such that $\ka(v+L)$ is an irrational winding in some torus in
$V/\Ga$).
\end{proof}

\begin{proof}[Proof of Theorem~\ref{reduc}]
If $\Ga$ is elliptic then ({\romannumeral1}) holds according to
Lemma~\ref{eqria} and known results on Euclidean space forms (see
\cite[Chapter 2]{Wolf}). Note that $\Ga$ acts in the spacelike
space $L^\bot$.

Let $\Ga$ be non-elliptic. Then $G$ has a nontrivial fixed point
$v_0\in\partial C$ by Lemma~\ref{fixed}, the space $W=v_0^\bot$ is
$G$-invariant and tangent to $\partial C$ at $v_0$, the
representation of $G$ in $V/W$ is trivial. Moreover, there is no
$G$-fixed point in $V\setminus W$ (otherwise, $G$ is elliptic
since each subspace that contains $v_0$ and is not contained in
$W$ intersects $\Int(C)$). Pick any $v_1\in\partial C$ such that
$\ell(v_0,v_1)=1$ and
set $L=\bbR v_0$, $N=v_1^\bot\cap W$ (thus $W,L,N$ are defined as
in ({\romannumeral2}) and (\ref{norml}) holds). By
Proposition~\ref{betow}, $\tau(\Ga)\subset W$. This makes possible
to apply Lemma~\ref{parab}. According to it, $\Ga$ acts in $W/L$
freely and discontinuously. Let $\phi$ be the quotient mapping.
Since $-\ell$ induces Euclidean structure in $W/L$, we may assume
that $\phi\circ\tau$ isomorphically embeds $\Ga$ to the vector
group $W/L$ (replacing $\Ga$ by its finite index subgroup if
necessary) as an uniform lattice in some linear subspace of $W/L$.
Due to Lemma~\ref{reafo}, the action of $\Ga$ is subject to
formulas of (\romannumeral2) (the vector $v_1$ may be removed);
Lemma~\ref{lastl} implies (\ref{freea}).

Conversely, let $\Ga$ be as in (\romannumeral 2) and (\ref{freea})
be true (the case (\romannumeral 1) is clear). Using relations
$G\subset\SO(\ell)$, $\ell(v_0,x)=\ell(v_0,ax)=0$, by a
straightforward calculation we get for any $v,u\in V$ and
$x\in\Ga$:
\begin{eqnarray}
\ell(\ga_x(u)-v,\,\ga_x(u)-v)
=\ell(x+l_0(v)ax,x+l_0(u)ax)+\dots,\label{negax}
\end{eqnarray}
where dots denote summands that are linear or constant on $x$. It
follows from definitions that $\ell$ is negative definite on $N$.
Due to (\ref{freea}), for each $v\in V$ the quadratic form
$\ell(x+l_0(v)ax,x+l_0(v)ax)$ is negative definite on $T$; if $u$
is sufficiently close to $v$ then the same is true for the form
\begin{eqnarray*}
q_{u,v}(x)=\ell(x+l_0(v)ax,\,x+l_0(u)ax).
\end{eqnarray*}
Therefore, there exist a neighbourhood $U$ of $v$ in $V$ and
$n_v\in\bbN$ such that for all $u\in U$ the number of $x$
satisfying the inequality
\begin{eqnarray*}
\ell(\ga_x(u)-v,\ga_x(u)-v)\geq0
\end{eqnarray*}
is less than $n_v$. Hence $\ka((v+C)\cap U)$ and $\ka((v-C)\cap
U)$ are closed in $\ka(U)\subset V/\Ga$. Thus, $V/\Ga$ is strictly
causal.

Let $S$ be the union of affine subspaces
\begin{eqnarray*}
S_t=(({\mathbf1}+ta)T)^\bot\cap X_t,\qquad t\in\bbR,
\end{eqnarray*}
where $X_t=\{v\in V:\,l_0(v)=t\}$.  Due to (\ref{freea}), $S$ is
homeomorphic to the vector space $V/T$. Let $v\in V_t$ and denote
by $u$ the orthogonal projection of $v$ to $({\mathbf1}+ta)T$ (the
definition is correct since $\ell$ is nondegenerate in $N$). There
exists the unique $x\in T$ such that $u+x+tax=0$. Since
$v_0,v_1\perp N$ and $\la(x)\in\SO(\ell)$, by (\ref{forla}) and
(\ref{forta}) this implies that each $T$-orbit in $V$ has
precisely one common point with $S$. Therefore, $V$ can be
realized as the trivial vector bundle over $S$ with the fibre
$T$. Thus $M$ is homeomorphic to $V/T\times T/\Ga$.
\end{proof}
\begin{proof}[Proof of Theorem~\ref{fourc}]
Let $M=V/\Ga$ be a 4-manifold of the type (\romannumeral2).
Suppose that $a\neq0$. Due to (\ref{freea}), this implies $N\neq
T$. Hence $\dim T=1$, $\dim N=2$, and $\Ga\cong\bbZ$. In the
coordinate system whose origin is removed to a point $\td o\in V$
the translation part of an affine transformation
$\ga(v)=\la(\ga)v+\tau(v)$ is replaced by
$$\td\tau(\ga)=\ga(\td o)-\td o=\tau(\ga)+\la(\ga)\td o-\td o$$
and $\la(\ga)$ does not change. A simple calculation with
(\ref{nuvot}), (\ref{lanua}), (\ref{tanua}) shows that there is a
point $\td o\in\bbR v_1$ such that $\td\tau(x)$ contains no
summands that are quadratic on $x\in\Ga$ (note that $\dim T=1$,
$\ell(ax,ax)\neq0$ if $x\neq0$, and that only coefficients at
$v_0$ have quadratic summands). In the proof of
Theorem~\ref{reduc} we did not specify the position of origin.
Hence we get the same formulas (with other $a$ and $v_1$ in
general) removing it to $\td o$. Then $\tau(x)$ is linear on $x$
and (\ref{tanua}) implies that $\ell(ax,x)=0$ and $\tau(x)=x$. Let
$e_1\in T$ generate $\Ga$. Multiplying $\ell$ by a suitable
positive number, we may assume $\ell(e_1,e_1)=-1$. There exists
the unique $t>0$ such that
$$\ell(\la(e_1)(tv_1)-tv_1,\,\la(e_1)(tv_1)-tv_1)=t^2\ell(ae_1,ae_1)=-1.$$
Set $e_3=tv_1$, $e_0=v_0/t$, $e_2=tae_1$. In this base we get
formulas of Example~\ref{newun} (note that $a$ must be replaced by
$ta$). Hence $M$ is homothetic to the manifold of this example.

Let $M'=V/\Ga'$ satisfy assumptions of the theorem. Replacing it
by a homothetic manifold, we may assume that $M'$ it is finitely
covered by the manifold $M=V/\Ga$ of Example~\ref{newun}. Then
$\Ga$ is a subgroup of finite index in $\Ga'$. By
Proposition~\ref{betow}, Lemma~\ref{parab} and
Corollary~\ref{iseuc}, $\Ga'$ acts freely and properly in
Euclidean 2-plane $W/L$ as a group of isometric transformations.
If $\Ga'$ acts by translations then we may apply
Lemma~\ref{reafo}. Then we have the setting above. Note that $\dim
T=1$ and $\Ga'\cong\bbZ$ in this case. Any other isometric
transformation that $\Ga'$ could contain is the composition of a
shift with the reflection with respect to the line of the shift.
Since we may apply the proven assertion to the subgroup of $\Ga'$
consisting of all translations in it, $\Ga'$ is generated by an
element $\ga$ of this type. Thus we may assume that $\Ga$ is the
group of Example~\ref{newun} and has index 2 in $\Ga'$; moreover,
that $\ga\circ\ga=\ga_1$, where $\ga_1$ is defined by
(\ref{defex}), (\ref{defon}). Set $g=\la(\ga)$ and $b=\tau(\ga)$.
Then
\begin{eqnarray*}
2b=e_1,\quad ge_0=e_0,\quad ge_1=e_1+se_0,\quad
ge_2=-e_2+re_0,\quad s,r\in\bbR.
\end{eqnarray*}
Hence $g^2e_2=e_2$; since $g^2=\la(1)$ and $\la(1)e_2=e_2+e_0$, we
get a contradiction.
\end{proof}

\begin{proof}[Proof of Proposition~\ref{disti}]
Extending the identical affine mapping to the universal covering
space we get a $G$-invariant positive definite quadratic form.
Hence the closure of $G$ is compact and $M$ is elliptic by
Lemma~\ref{eleco}. Conversely, if $M$ is elliptic then $G$ admits
a $G$-invariant positive definite form that defines Euclidean
structure in $M$. Therefore, $M$ is elliptic if and only if 1)
holds.

If $M$ is not elliptic then it can be finitely covered by $V/\Ga$
with $\Ga$ as in (\romannumeral2), where $a\neq0$.  Let $x\in\Ga$
and $ax\neq0$. Then there exists $c>0$ such  that
\begin{eqnarray}\label{lessz}
u,v\in V,\ -l_0(u)>c,\ l_0(v)>c\quad\Longrightarrow\quad
\ell(x+l_0(v)ax,\,x+l_0(u)ax)>0.
\end{eqnarray}
Due to~(\ref{negax}), for $u,v$ as above and all sufficiently
large $n\in\bbZ$
\begin{eqnarray*}
\ell(\ga_{nx}(u)-v,\,\ga_{nx}(u)-v)>0.
\end{eqnarray*}
Fixing $v$, we get that the past of $\ka(v)$ contains all points
$\ka(u)$ with $l_0(u)<-c$. This set is invariant with respect to
the translations along $W$, in particular, along $L$. Hence the
past of $\ka(v)$ contains lightlike straight lines. Thus $M$ is
elliptic if 2) is true.  Similar arguments prove 3).

Let $\Ga$ be of the type (\romannumeral2) with $a\neq0$, $M$ be
finitely covered by $V/\Ga$ and suppose that $M$ admits the Cauchy
hypersurface $S$. Then for any $v\in V$ the line $\ka(v+L)$ has a
common point with $S$. Hence  we may assume that $p=\ka(v)\in S$
and $l_0(v)>c$, where $c$ is as in (\ref{lessz}). Then $P_p$
contains all lines $\ka(u+L)$, where $u\in V$ satisfies
$l_0(u)<-c$.  These lines cannot intersect $S$ since $S$ separates
$P_p$ and $F_p$ being a Cauchy hypersurface, and we get a
contradiction. Thus 4) implies that $M$ is elliptic.

Since $\codim T\geq2$ in (\romannumeral2) and $M$ is homeomorphic
to $V/T\times T/\Ga$, it cannot have two ends.

Let $M$ be elliptic. Then $G=\la(\Ga)$ is finite and keeps $C$.
Hence there exists $G$-fixed vector $v_0\in\Int(C)$ and
$\Ga$-invariant spacelike hyperplane $H=v_0^\bot$ in $V$.
Evidently, $\ka(H)$ is a Cauchy hypersurface. Hence $M$ satisfies
4). The function $\ell(v_0,v)$ is $\Ga$-invariant and strictly
increasing along all lightlike lines; moreover, it is not bounded
on each of them. Thus 2) and 3) are true if $M$ is elliptic.
\end{proof}

\section*{Acknowledgements}
The initial version of the paper did not contain
Example~\ref{newun}, Theorem~\ref{fourc} 
and a part of Theorem~\ref{reduc}. We thank the referee for useful
comments that stimulate more detailed study of strictly causal
manifolds and for making us aware of the paper \cite{Barbot}.


\begin{thebibliography}{9999}
\bibitem{Abels}
Abels, H.: Properly Discontinuous Groups of Affine Transformations: A Survey,
{\it Geometriae Dedicata} 87 (2001), 309-–333.
\bibitem{Barbot}
Barbot, T.: Globally hyperbolic flat space–times, {\it Journal of
Geometry and Physics} 53 (2005), 123–-165.
\bibitem{Beem} Beem, J.K., Ehrlich, P.E.: Global Lorentzian geometry,
{\it Monographs and Textbooks in Pure and Applied Mathematics},
2nd ed., vol. 202, Marcel Dekker, New York, 1996.
\bibitem{CDGM03}
Charette, V., Drumm, T., Goldman, W. and Morill M.:
Complete Flat Affine and Lorentzian Manifolds, {\it Geometriae Dedicata} 97 (2003), 187-–198.
\bibitem{HoEl} Ellis, G., Hawking, S.: The large scale structure of space–time,
{\it Cambridge Monographs on Mathematical Physics, No. 1,}
Cambridge University Press, London, New York, 1973.
\bibitem{Fried}
Fried, D.: Flat Spacetimes, {\it J. Diff. Geom.} 26 (1987), 385--396.
\bibitem{Mar83} Margulis, G.: Free properly discontinuous groups of affine transformations,
{\it Dokl. Akad. Nauk SSSR} 272 (1983), 937--940.
\bibitem{Mar87} Margulis, G.: Complete affine locally flat manifolds with a free fundamental
group, {\it J. Soviet Math.} 134(1987), 129--134.
\bibitem{Milnor} Milnor, J.: On fundamental groups of complete affinely flat manifolds,
{\it Adv. in Math.} 25(1977), 178--187.
\bibitem{Wolf} Wolf, J.A: Spaces of constant curvature, University of Californnia, Berkley,
California, 1972.



\end{thebibliography}
\end{document}